\documentclass{amsart}
\usepackage{hyperref}
\usepackage{amsmath,amsthm,amscd,amssymb}
\usepackage{enumerate}
\usepackage{curves}

\newcommand*{\mailto}[1]{\href{mailto:#1}{\nolinkurl{#1}}}

\newtheorem{theorem}{Theorem}[section]
\newtheorem{lemma}[theorem]{Lemma}
\newtheorem{corollary}[theorem]{Corollary}

\newtheorem{hypothesis}[theorem]{Hypothesis {\bf H.}\hspace*{-0.6ex}}

\numberwithin{equation}{section}

\newcommand{\C}{\mathbb{C}}
\newcommand{\R}{\mathbb{R}}
\newcommand{\N}{\mathbb{N}}
\newcommand{\E}{\mathrm{e}}
\newcommand{\I}{\mathrm{i}}

\newcommand{\Res}{\mathop{\rm Res}}
\newcommand{\sign}{\mathop{\rm sign}}

\renewcommand{\Im}{\mathop{\rm Im}}

\newcommand{\beq}{\begin{equation}}
\newcommand{\eeq}{\end{equation}}
\newcommand{\bal}{\begin{align}}
\newcommand{\eal}{\end{align}}
\newcommand{\nn}{\nonumber}

\newcommand{\la}{\lambda}

\def\XXint#1#2#3{{\setbox0=\hbox{$#1{#2#3}{\int}$}
     \vcenter{\hbox{$#2#3$}}\kern-.5\wd0}}

\numberwithin{equation}{section}

\begin{document}

\title[The transformation operator for Schr\"odinger operators]{The transformation operator for Schr\"odinger operators on almost periodic infinite-gap backgrounds}

\author[K. Grunert]{Katrin Grunert}
\address{Faculty of Mathematics\\ Nordbergstrasse 15\\ 1090 Wien\\ Austria}
\email{\mailto{katrin.grunert@univie.ac.at}}
\urladdr{\url{http://www.mat.univie.ac.at/~grunert/}}

\thanks{Research supported by the Austrian Science Fund (FWF) under Grant No.\ Y330.}
\thanks{J. Diff. Eq. {\bf 250}, 3534--3558 (2011)}


\begin{abstract}
We investigate the kernels of the transformation operators for one-dimensional Schr\"odinger operators with potentials, which are asymptotically close to Bohr almost periodic infinite-gap potentials.
\end{abstract}
  
  \maketitle

\section{Introduction}

In the investigation of direct and inverse scattering problems transformation operators play an important role. They first showed up in the context of generalized shift operators in the work of Delsarte \cite{D}. They have also been investigated by Levitan  \cite{L4} and were constructed for arbitrary Sturm-Liouville equations by Povzner \cite{P}. Afterwards transformation operators have been applied for the first time when considering inverse spectral problems, for example by Marchenko \cite{M}, who noticed that the spectral function of a Sturm-Liouville operator determines the operator uniquely in \cite{M2}. Soon after that Gel'fand and Levitan \cite{GL} found a method of recovering Sturm-Liouville equations from its spectral functions, using the transformation operator techniques.\\
Another important step was the introduction of transformation operators, which preserve the asymptotic behavior of solutions at infinity by Levin \cite{Le}. Since that these transformation operators are the main tool for solving different kinds of scattering problems, mainly in the case of constant backgrounds. They have been studied for periodic infinite-gap backgrounds by Firsova \cite{F2}, \cite{F1} and they have been recently investigated in the finite-gap case by Boutet de Monvel, Egorova, and Teschl \cite{BET}. \\
 In the present work we propose a complete investigation of the transformation operators for Bohr almost periodic infinite-gap backgrounds, which belong to the so-called Levitan class. It should be noted, that this class includes as a special case the set of smooth, periodic infinite-gap operators.  

 To set the stage, we need

\begin{hypothesis}\label{Hyp}
Let
\begin{equation}\nn 
0\leq E_0< E_1<  \dots  <E_n< \dots
\end{equation}
be an increasing  sequence of points on the real axis which satisfies the following conditions:
\begin{enumerate}
\item
for a certain $l>1$, $\sum_{n=1}^\infty (E_{2n-1})^l (E_{2n}-E_{2n-1}) <\infty$ and 
\item 
$E_{2n+1}-E_{2n-1} > C n^{\alpha}$, where $C$ and $\alpha$ are some fixed, positive constants.
\end{enumerate}
\end{hypothesis}

We will call, in what follows, the intervals $(E_{2j-1},E_{2j})$ for $j=1,2,\dots$ gaps. In each closed gap $[E_{2j-1},E_{2j}]$ , $j=1,2,\dots$, we choose a point $\mu_j$ and an arbitrary sign $\sigma_j\in\{\pm 1\}$.

Consider next the system of differential equations for the functions $\mu_j(x)$, $\sigma_j(x)$, $j=1,2,...$, which is an infinite analogue of the  well-known Dubrovin equations, given by 
\begin{align}\label{repmujx}
\frac{d \mu_j(x)}{d x}=
& -2\sigma_j(x)\sqrt{-(\mu_j(x)-E_0)}\sqrt{\mu_j(x)-E_{2j-1}}\sqrt{\mu_j(x)-E_{2j}} \\ \nn 
& \times\prod_{k=1, k\not =j}^\infty  \frac{\sqrt{\mu_j(x)-E_{2k-1}}\sqrt{\mu_j(x)-E_{2k}}}{\mu_j(x)-\mu_k(x)}
\end{align}
with initial conditions $\mu_j(0)=\mu_j$ and $\sigma_j(0)=\sigma_j$, $j=1,2,\dots$
\footnote{We will use the standard branch cut of the square root in the domain $\mathbb C\setminus\mathbb R_+$ with $\Im\sqrt z>0$.}. Levitan \cite{L2}, \cite{L3}, and \cite{L}, proved, that
 this system of differential equations is uniquely solvable, that the solutions $\mu_j(x)$, $j=1,2,\dots$ are continuously differentiable and satisfy $\mu_j(x)\in [E_{2j-1}, E_{2j}]$ for all $x\in\R$. Moreover, these functions $\mu_j(x)$, $j=1,2,\dots$ are Bohr almost periodic
\footnote{
For informations about almost periodic functions we refer to \cite{LZ}.}.
Using the trace formula (see for example \cite{L})
\beq\label{defp}
p(x)=E_0+\sum_{j=1}^\infty (E_{2j-1}+E_{2j}-2\mu_j(x)),
\eeq
we see that also $p(x)$ is real Bohr almost periodic. The operator 
\begin{equation}\label{deftiL}
 \tilde L:= -\frac{d^2}{dx^2}+p(x), \quad \text{dom}(\tilde L)=H^{2}(\R),
\end{equation}
in $L^2(\R)$, is then called an almost periodic infinite-gap Schr\"odinger operator of the Levitan class. The spectrum of $\tilde L$ is purely absolutely continuous and of the form
\begin{equation}\nn
\sigma=[E_0,E_1]\cup \dots\cup [E_{2j},E_{2j+1}] \cup \dots,
\end{equation}
and has spectral properties analogous to the quasi-periodic finite-gap Schr\"odinger operator.
In particular, it is completely defined by the series $\sum_{j=1}^\infty (\mu_j,\sigma_j)$, which we call the Dirichlet divisor. Analogously to the finite--gap case this divisor is connected to a Riemann surface of infinite genus, which is associated to the function $Y^{1/2}(z)$, where  
\beq\label{defY}
Y(z)=-(z-E_0)\prod_{j=1}^\infty  \frac{(z-E_{2j-1})}{E_{2j-1}}\frac{(z-E_{2j})}{E_{2j-1}}, 
\eeq
and where the branch cuts are taken along the spectrum. 
It is known, that the Schr\"odinger equation 
\begin{equation}\nn
\Big(-\frac{d^2}{dx^2}+p(x)\Big)y(x)=z y(x)
\end{equation}
with any continuous, bounded potential $p(x)$ has two Weyl solutions $\psi_\pm(z,x)$ normalized by 
\begin{equation*}
\psi_\pm(z,0)=1 \quad \text{ and } \quad \psi_\pm(z,.)\in L^2(\R_\pm),\quad \text{ for } z\in\C\backslash\sigma.
\end{equation*}
In our case of Bohr almost periodic potentials of the Levitan class, these solutions have complementary properties similar to properties of the Baker-Akhiezer functions in the finite-gap case. We will briefly discuss them in the next section. 

The objects of interest, for us, are the Jost solutions of the one-dimensional Schr\"odinger operator $L$ in $L^2(\R)$
\beq\label{defL}
L:=-\frac{d^2}{dx^2}+q(x), \quad \text{dom}(L)=H^{2}(\R), 
\eeq
with the real-valued potential $q(x)\in C(\R)$ satisfying the following condition
\beq\label{moment}
\int_\R (1+\vert x\vert^2 ) \vert q(x)-p(x)\vert dx<\infty. 
\eeq

We will prove the following result.
\begin{theorem}\label{mainth}
 Assume Hypothesis~\ref{Hyp}. Let $p(x)$, defined as in \eqref{defp}, belong to the Levitan class, and $q(x)$ satisfy \eqref{moment}, then the Jost solutions $\phi_\pm(z,x)$ can be represented in the following form 
\begin{equation*}
\phi_\pm(z,x)=\psi_\pm(z,x)\pm\int_x^{\pm\infty} K_\pm(x,y)\psi_\pm(z,y)dy,
\end{equation*}
where the solutions of \eqref{inteq:K}, are real valued, continuously differentiable with respect to both parameters and for $\pm y>\pm x$ they satisfy 
\begin{equation*}
\vert K_\pm(x,y)\vert \leq \pm C_\pm(x)\int_{\frac{x+y}{2}}^{\pm\infty}\vert q(s)-p(s)\vert ds.
\end{equation*}
Here $C_\pm(x)$ are continuous positive functions, which are monotonically decreasing as $x\to\pm\infty$.  
\end{theorem}

Here it should be pointed out that a subset of the operators belonging to the Levitan class consists of operators with periodic potentials. Assume that the sequence $\{E_j\}_{j=1}^\infty$ fulfills Hypothesis~\ref{Hyp}, then there is a criterion when this sequence is the set of band edges of the spectrum of some Schr\"odinger operator with periodic potential $p(x+a)=p(x)\geq 0$ with $p\in W_2^3[0,a]$. Namely, Marchenko and Ostrovskii proved in \cite{MO}, that 
\begin{equation*}
 p\in W_2^k[0,a] \quad \text{ iff } \quad \sum_{j=1}^\infty j^{2k+2}(\sqrt{E_{2j}-E_0}-\sqrt{E_{2j-1}-E_0})^2<\infty,
\end{equation*}
for $k=0,1,\dots$. 
As it is well known that in the periodic case $E_{2j-1}=j^2+O(1)$ and $E_{2j}=j^2+O(1)$ as $j\to\infty$, we obtain for large $j$ that 
\begin{align}\nn
 E_{2j-1}(E_{2j}-E_{2j-1})
& \leq E_{2j-1}(\sqrt{E_{2j}-E_0}+\sqrt{E_{2j-1}-E_0})\times\\ \nn
& \quad (\sqrt{E_{2j}-E_0}-\sqrt{E_{2j-1}-E_0})\\ \nn
& \leq 2j^3(\sqrt{E_{2j}-E_0}-\sqrt{E_{2j-1}-E_0})=:I_j.  
\end{align}
Since Hypothesis~\ref{Hyp} is satisfied, we have $\sum_{j=1}^\infty j^{2k-4}I_j^2<\infty$ for $k>2$ and hence the Cauchy inequality implies that $\sum_{j=1}^\infty I_j<\infty$ in this case. This means, that Hypothesis~\ref{Hyp} is satisfied for any $a$-periodic potential $p(x)\geq 0$ with $p\in W_2^k[0,a]$ for $k>2$.

The Jost solutions $\phi_\pm(z,x)$ are bounded for $z\in \sigma(\tilde L)$ and in particular there is no subordinate solution for $z$ in the interior of $\sigma(\tilde L)$. Hence the essential spectrum is purely absolutely continuous (cf.\ \cite{GMZ}, \cite{GZ}, \cite{G1},\cite{G2}, \cite{GP}, and \cite[Thm.~9.31]{Te}):

\begin{corollary}
The essential spectrum of $L$ is purely absolutely continuous
\begin{equation}
\sigma_{ac}(L) = \sigma(\tilde L), \qquad \sigma_{sc}(L) = \emptyset.
\end{equation}
The point spectrum is confined to the closed gaps $\overline{\R\backslash \sigma(\tilde L)}$.
\end{corollary}

Criteria that there are only finitely many in each closed gap can be found in \cite{K}, \cite{KT}, \cite[Thm.~6.12]{RK}.

Another application is given by the scattering problem for the operators \eqref{deftiL} and \eqref{defL}, which has not been considered so far. Such scattering problems arise in various physical applications, for example, when studying properties of the alloy of two different semi-infinite one-dimensional crystals. A more detailed discussion of the history of such problems and further references to the literature can be found in \cite{GNP}. Rather recently the scattering problem has been solved in the case of Schr\"odinger operators with steplike finite-gap potentials in \cite{BET} and these results have been applied to solve the Korteweg--de Vries equation in \cite{EGT} and \cite{ET} (see also \cite{ET2,ET3,MT}). Even in the case of periodic backgrounds scattering theory has been developed \cite{F1} and the time evolution of the scattering data for the Korteweg--de Vries equation has been computed \cite{F3}. All these investigations use transformation operators as a main tool.

\section{Background Schr\"odinger operators}

In this section we want to summarize some facts for the background Schr\"odinger operator of Levitan class. 
We present these results, obtained in \cite{L}, \cite{SY}, and \cite{SY2}, in a form similar to the finite-gap case used in \cite{BET} and \cite{GH}.
 
Let $\tilde L$ be as in \eqref{deftiL}.
Denote by 
$s(z,x)$, $c(z,x)$ the sine- and cosine-type solutions of the corresponding equation 
\beq\label{eq:schr}
\left( -\frac{d^2}{dx^2} +p(x) \right) y(x)=zy(x), \quad z\in\C,
\eeq
associated with the initial conditions 
\begin{equation*}
s(z,0)=c^\prime(z,0)=0, \quad c(z,0)=s^\prime(z,0)=1,
\end{equation*}
where prime denotes the derivative with respect to $x$.
Then $c(z,x)$, $c^\prime(z,x)$, $s(z,x)$, and $s^\prime(z,x)$ are entire with respect to $z$. 
They can be represented in the following form
\begin{equation*}
c(z,x)=\cos(\sqrt{z}x)+\int_0^x\frac{\sin(\sqrt{z}(x-y))}{\sqrt{z}}p(y)c(z,y)dy,
\end{equation*}
\begin{equation*}
s(z,x)=\frac{\sin(\sqrt{z}x)}{\sqrt{z}}+\int_0^x \frac{\sin(\sqrt{z}(x-y))}{\sqrt{z}}p(y)s(z,y)dy.
\end{equation*}

The background Weyl solutions are given by 
\begin{equation} \label{1.5}
\psi_{\pm}(z,x) = c(z,x) + m_{\pm}(z,0) s(z,x),
\end{equation}
where 
\beq\label{mfunc}
m_{\pm}(z,x)=\frac{\psi_\pm^\prime (z,x)}{\psi_\pm(z,x)}=\frac{H(z,x)\pm  Y^{1/2}(z)}{G(z,x)},
\eeq 
are the Weyl functions of $\tilde L$ (cf. \cite{L}),
where $Y(z)$ is defined by (\ref{defY}), 
\beq\label{repH}
G(z,x)=\prod_{j=1}^\infty\frac{z-\mu_j(x)}{E_{2j-1}}, \quad \text{ and } \quad  
H(z,x)=\frac{1}{2}\frac{d}{d x}G(z,x).
\eeq 
Using (\ref{repmujx}) and (\ref{repH}), we have 
\begin{equation}\label{repHG}
H(z,x)=\frac{1}{2}\frac{d}{d x}G(z,x)= G(z,x) \sum_{j=1}^\infty \frac{\sigma_j(x) Y^{1/2}(\mu_j(x))}{\frac{d}{dz} G(\mu_j(x),x)(z-\mu_j(x))}.
\end{equation}
The Weyl functions $m_\pm(z,x)$ are Bohr almost periodic as the following argument shows:
For each $j\in \N$ the functions $\mu_j(x)$ are almost periodic and hence, as a finite product of almost periodic functions is again almost periodic, also 
\beq
G_n(z,x)=\prod_{j=1}^n \frac{z-\mu_j(x)}{E_{2j-1}}
\eeq
is almost periodic for fixed $z\in \C$.
Moreover, we have 
\beq\label{quasiG}
\vert G_n(z,x)-G(z,x) \vert 
=\left\vert \prod_{j=1}^n \frac{z-\mu_j(x)}{E_{2j-1}}\Big(1 -\prod_{j=n+1}^\infty \frac{z-\mu_j(x)}{E_{2j-1}}\Big)\right\vert,
\eeq
 and for every fixed $z\in \C$ there exists a $m\in\N$ such that $\vert z\vert \leq E_{2m-1}$.  Then for $n>m$, we obtain on the one hand 
\begin{equation} 
\exp \Big(-\sum_{j=n+1}^\infty \frac{\vert z\vert }{E_{2j-1}-\vert z \vert }\Big)
\leq \prod_{j=n+1}^\infty \frac{\mu_j(x)-\vert z \vert }{E_{2j-1}}
 \leq \prod_{j=n+1}^\infty \vert \frac{z-\mu_j(x)}{E_{2j-1}}\vert, 
\end{equation}
and on the other hand 
\begin{align} \nn 
\prod_{j=n+1}^\infty \vert \frac{z-\mu_j(x)}{E_{2j-1}}\vert 
& \leq \prod_{j=n+1}^\infty \frac{\mu_j(x)+\vert z\vert }{E_{2j-1}}\\
&\leq \exp \Big(\frac{1}{E_{2n+1} }\sum_{j=n+1}^\infty (E_{2j}-E_{2j-1})+\vert z\vert \sum_{j=n+1}^\infty \frac{1}{E_{2j-1}}\Big),
\end{align}
where we used that $\log(1+x)\leq x$ and $\log(1-x)\geq \frac{-x}{1-x}$ for $x>0$. Noticing that the second condition in Hypothesis~\ref{Hyp} implies that $\sum_{j=1}^\infty \frac{1}{E_{2j-1}}$ converges, all terms are well defined, and it follows that the product $\prod_{j=n+1}^\infty \frac{z-\mu_j(x)}{E_{2j-1}}$ converges to $1$ as $n\to\infty$. Furthermore, $\prod_{j=1}^\infty \frac{z-\mu_j(x)}{E_{2j-1}}$ is uniformly bounded with respect to $x$ for any fixed $z$. 
As all our estimates are independent of $x$, we have that $G_n(z,x)$ converges uniformly for fixed $z$ against $G(z,x)$ and thus, the function $G(z,x)$ is almost periodic with respect to $x$.  
Furthermore by definition $\frac{H(z,x)}{G(z,x)}=\frac{1}{2}\frac{G^\prime(z,x)}{G(z,x)}=\frac{1}{2}(\log(G(z,x)))^\prime$ and therefore $\frac{H(z,x)}{G(z,x)}$ is also almost periodic, where we use that $\log (G(z,x)) \not = 0$ for $z\not\in [E_1,E_2]\cup \dots\cup [E_{2j-1},E_{2j}]\cup\dots $ together with \cite[Property 3,4,5]{LZ}, and hence $m_\pm(z,x)$ are also almost periodic functions.

\begin{lemma}\label{lempsi}
The background Weyl solutions, for $z\in\C$, can be represented in the following form 
\begin{align}\label{intrep}
\psi_{\pm}(z,x) =\exp\left(\int_0^x m_{\pm}(z,y)dy\right)= \left( \frac{G(z,x)}{G(z,0)}\right)^{1/2}\exp\left( \pm \int_0^x \frac{Y^{1/2}(z)}{G(z,y)}dy \right).
\end{align}
If for some $\varepsilon>0$, $\vert z-\mu_j(x)\vert > \varepsilon$ for all $j\in\N$ and $x\in\R$, then the following holds:
For any $1>\delta>0$ there exists an $R>0$ such that 
\begin{equation}\label{psiexp}
 \vert\psi_{\pm}(z,x)\vert \leq \E^{\mp(1-\delta )x\Im(\sqrt{z})}\Big(1+\frac{D_R}{\vert z\vert}\Big), \text{ for any } \vert z\vert \geq R, \phantom{o}\pm x>0.   
\end{equation}
where $D_R$ denotes some constant dependent on $R$. 
\end{lemma}

\begin{proof}
 First we will show that 
\beq\label{fpm}
f_\pm(z,x)=\left( \frac{G(z,x)}{G(z,0)}\right)^{1/2}\exp\left( \pm \int_0^x \frac{Y^{1/2}(z)}{G(z,y)}dy \right)
\eeq
fulfills 
\beq\label{diff f}
\left(-\frac{d^2}{dx^2}+p(x)\right)y(x)=zy(x).
\eeq
Using \eqref{mfunc} we obtain that $f^\prime_\pm=m_\pm f_\pm$ and $f^{\prime\prime}_\pm = (m^\prime_\pm+m^2_\pm)f_\pm$. Hence 
\eqref{diff f} will be satisfied if and only if 
\begin{equation}
 m^\prime_\pm+m^2_\pm =p(x)-z. 
\end{equation}
This can be checked using the following relations, which are proved in \cite{L},
\begin{equation}\label{eq:GNHY}
 G(z,x)N(z,x)+H(z,x)^2=Y(z), 
\end{equation}
where 
\begin{equation}
 N(z,x)=-(z-\tau_0(x))\prod_{j=1}^\infty \frac{z-\tau_j(x)}{E_{2j-1}},
\end{equation}
with $\tau_0(x)\in (-\infty,E_0]$ and $\tau_j(x)\in[E_{2j-1},E_{2j}]$
and 
\begin{equation}
 \frac{d^2}{dx^2}G(z,x)= 2((p(x)-z)G(z,x)-N(z,x)). 
\end{equation}
Moreover, for $z$ outside an $\varepsilon$ neighborhood of the gaps we have $f(z,0)=1$, and we conclude that 
\beq
\frac{G(z,x)}{G(z,0)} =\prod_{j=1}^\infty \frac{z-\mu_j(x)}{z-\mu_j(0)}
 = \exp\Big(\sum_{j=1}^\infty \log\Big(1+\frac{\mu_j(0)-\mu_j(x)}{z-\mu_j(0)}\Big)\Big).
\eeq
Thus we obtain 
\beq
\left\vert\frac{G(z,x)}{G(z,0)}\right\vert \leq \exp\Big(\sum_{j=1}^\infty \log \Big(1+\left\vert \frac{\mu_j(0)-\mu_j(x)}{z-\mu_j(0)}\right\vert\Big)\Big)
 \leq \exp\Big(\sum_{j=1}^\infty \left\vert \frac{\mu_j(0)-\mu_j(x)}{z-\mu_j(0)}\right\vert\Big),
\eeq
where we used that $\log(1+x)\leq x$ for $x>0$. 
Moreover
\beq
\left\vert \frac{\mu_j(0)-\mu_j(x)}{z-\mu_j(0)}\right\vert 
=\frac{1}{\left\vert \frac{z-\mu_j(0)}{\mu_j(0)-\mu_j(x)}\right\vert},
\eeq
which implies for $\vert z \vert \leq 2E_{2j}$ that 
\beq
\left\vert \frac{(\mu_j(0)-\mu_j(x))z}{z-\mu_j(0)}\right\vert \leq \left\vert \frac{2(\mu_j(0)-\mu_j(x))E_{2j}}{\varepsilon}\right\vert.
\eeq
For $\vert z \vert > 2E_{2j}$ we can estimate the terms by 
\begin{align} 
\left\vert \frac{(\mu_j(0)-\mu_j(x))z}{z-\mu_j(0)}\right \vert & 
\leq \vert \mu_j(0)-\mu_j(x)\vert \left\vert\frac{1}{1-\frac{\mu_j(0)}{z}}\right\vert 
\leq \vert \mu_j(0)-\mu_j(x) \vert\frac{1}{ 1-\left\vert\frac{\mu_j(0)}{z}\right\vert }\\ \nn
& \leq \vert \mu_j(0)-\mu_j(x)\vert \frac{1}{1-\left\vert \frac{E_{2j}}{z}\right\vert}
\leq 2\vert \mu_j(0)-\mu_j(x)\vert.
\end{align}
Combining the estimates from above, we obtain
\begin{align}
 \left\vert \frac{(\mu_j(0)-\mu_j(x))z}{z-\mu_j(0)}\right\vert 
&\leq 2\max (1,\frac{E_{2j}}{\varepsilon})\vert \mu_j(0)-\mu_j(x)\vert\\ \nn
& \leq 2\max (1,\frac{E_{2j}}{\varepsilon})(E_{2j}-E_{2j-1}),
\end{align}
and for any fixed $\varepsilon>0$ there exists a $k$ independent of $x$ and $z$ such that $\frac{E_{2n}}{\varepsilon}>1$ for all $n>k$, and therefore
\beq
\left\vert \frac{G(z,x)}{G(z,0)}\right\vert 
\leq \exp\Big(\frac{1}{\vert z \vert}\sum_{j=1}^\infty\left\vert \frac{(\mu_j(0)-\mu_j(x))z}{z-\mu_j(0)}\right\vert\Big)
\leq \exp\Big(C\frac{1}{\vert z \vert}\Big).
\eeq
where $C$ is a constant independent of $x$ and $z$. 
Analogously one can now investigate $\frac{Y^{1/2}(z)}{G(z,x)}$. 
Using 
\begin{equation}
 Y^{1/2}(z)=\I \sqrt{z-E_0}\prod_{j=1}^\infty\frac{\sqrt{z-E_{2j-1}}\sqrt{z-E_{2j}}}{E_{2j-1}},
\end{equation}
where the roots are defined as follows
\begin{equation}
 \sqrt{z-E}=\sqrt{\vert z-E\vert }\E^{\I\arg (z-E)/2},
\end{equation}
together with
$\frac{Y^{1/2}(z)}{G(z,x)}$ is a Herglotz function and 
\begin{equation}
 \sqrt{z-E_0}=\sqrt{z}(1+O(\frac{1}{z})), \quad \text{ as } z\to\infty, 
\end{equation}
we obtain the following estimate, which is uniform with respect to $x$,
\begin{equation}
 \frac{Y^{1/2}(z)}{G(z,x)}=\I\sqrt{z}(1+O(\frac{1}{z})), \quad \text{ as } z\to\infty. 
\end{equation} 
Using now that $\int_0^x \frac{Y^{1/2}(z)}{G(z,\tau)}d\tau = x\frac{Y^{1/2}(z)}{G(z,\zeta)}$, where $\zeta\in (0,x)$ by the mean value theorem, we finally obtain that $f_\pm(z,x)$ has the following asymptotic expansion outside a small neighborhood of the gaps as $z\to\infty$
\begin{equation}\label{fexp}
 f_\pm(z,x)=\E^{\pm\I\sqrt{z}x(1+O(\frac{1}{z}))}\Big(1+O(\frac{1}{z})\Big),
\end{equation}
where the $O(\frac{1}{z})$ terms are uniformly bounded with respect to $x$ and 
where we use the branch cut of the square root in the domain $\C\backslash\R_+$ with $\Im (\sqrt{z})>0$. 
Thus $f_\pm(z,.) \in L^2(\R_\pm)$ for $z\in\C$ outside a small neighborhood of the gaps and  
therefore away from the gaps $f_\pm(z,x)$ must coincide with $\psi_\pm(z,x)$. This implies in particular that $f_\pm(z,x)$ is holomorphic as a function of $z$ in the domain $\C\backslash \{z\in\C \ |\ \Im (z)\leq \varepsilon\}$.
For any point $z_0\in \C$ with $\Im(z_0)=2\varepsilon$, we can replace $f_\pm(z,x)$ by its Taylor series around the point $z_0$, which has at least radius of convergence $\varepsilon$. This Taylor series must coincide with the Taylor expansion of $\psi_\pm(z,x)$ around $z_0$, which converges on the disc with radius $2\varepsilon$ as $\psi_\pm(z,x)$ are holomorphic on the domain $\C\backslash\R$. Thus $f_\pm(z,x)$ is holomorphic on the domain $\C\backslash\R$ and hence $\psi_\pm(z,x)=f_\pm(z,x)$ on $\C\backslash\R$. Moreover, this implies that $f_\pm(z,x)$ and $\psi_\pm(z,x)$ must coincide on the real axis and hence, we finally get that $\psi_\pm(z,x)$ can be represented for any $z\in\C$ by (\ref{intrep}) and \eqref{fexp} implies \eqref{psiexp}.  
\end{proof}

As the spectrum consists of infinitely many bands, let us cut the complex plane along the spectrum $\sigma$ and
denote the upper and lower sides of the cuts by $\sigma^u$ and
$\sigma^l$. The corresponding points on these cuts will be denoted by
$z^u$ and $z^l$, respectively. In particular, we introduce the notation
\[
f(z^u) := \lim_{\varepsilon\downarrow0} f(z+\I\varepsilon),
\qquad f(z^l) := \lim_{\varepsilon\downarrow0}
f(z-\I\varepsilon), \qquad z\in\sigma,
\]
whenever the limits exist.
Define the Green function (see e.g. \cite{BG}, \cite{C}, \cite{KK}, and \cite{SY}) 
\begin{equation}\label{1.88}
g(z)= -\frac{G(z,0)}{2
Y^{1/2}(z)},
\end{equation}
where the branch of the square root is chosen in such a way that
\begin{equation}\label{1.8}
\frac{1}{\I} g(z^u) = \Im(g(z^u))  >0 \quad
\mbox{for}\quad z\in\sigma,
\end{equation}
then we obtain after a short calculation 
\beq\label{1.62}
W(\psi_-(z), \psi_+(z))=m_+(z)-m_-(z)=- g(z)^{-1},
\eeq
where $W(f,g)(x)=f(x)g^\prime(x)-f^\prime(x)g(x)$ denotes the usual Wronskian determinant.

For every Dirichlet eigenvalue $\mu_j=\mu_j(0)$, the Weyl functions $m_{\pm}(z)$ might have singularities. 
If $\mu_j$ is in the interior of its gap, precisely one Weyl function $m_+$ or $m_-$
will have a simple pole. Otherwise, if $\mu_j$ sits at an edge, both will have
a square root singularity. Hence we divide the set of poles accordingly:
\begin{align*}
M_+ &=\{ \mu_j\mid\mu_j \in (E_{2j-1},E_{2j}) \text{ and } m_+ \text{ has a simple pole}\},\\
M_- &=\{ \mu_j\mid\mu_j \in (E_{2j-1},E_{2j}) \text{ and } m_- \text{ has a simple pole}\},\\
\hat M &=\{ \mu_j\mid\mu_j \in \{E_{2j-1},E_{2j}\} \}.
\end{align*}

In particular, the following properties of the Weyl solutions are valid  (see, e.g. \cite{CL}, \cite{L}, \cite{SY}, \cite{Te}):

\begin{lemma}\label{lem1.1}
The Weyl solutions have the following properties:
\begin{enumerate}[\rm(i)]
\item
The functions $\psi_{\pm}(z,x)$ are holomorphic
as a function of $z$ in the domain $\mathbb{C}\setminus (\sigma\cup M_{\pm})$,
real valued on the set $\mathbb{R}\setminus \sigma$, and have
simple poles at the points of the set $M_{\pm}$.
Moreover, they are  continuous up to the boundary $\sigma^u\cup \sigma^l$ except at the points in $\hat M$ and
\begin{equation}\label{1.10}
\psi_+(z^u) = \psi_-(z^l) =\overline{\psi_+(z^l)},
\quad z\in\sigma.
\end{equation}
For $E \in \hat M$ the Weyl solutions satisfy
\begin{equation}\label{behpsi}
\psi_{\pm}(z,x)=O\left(\frac{1}{\sqrt{z-E}}\right), \quad
\mbox{as } z\to E\in \hat M,
\end{equation}
where the $O((z-E)^{-1/2})$ term is independent of $x$.

\noindent
The same applies to $\psi'_{\pm}(z,x)$.
\item
The functions $\psi_{\pm}(z,x)$ form an orthonormal basis on the spectrum with respect to the weight 
\begin{equation}\label{1.12}
d\rho(z)=\frac{1}{2\pi\I}g(z) dz,
\end{equation}
and any $f(x)\in L^2(\R)$ can be expressed through
\beq\label{1.14}
f(x)=\oint_{\sigma} \left(\int_\R f(y)\psi_+(z,y)dy\right)\psi_-(z,x)d\rho(z).
\eeq
Here we use the notation 
\beq
\oint_{\sigma}f(z)d\rho(z) := \int_{\sigma^u} f(z)d\rho(z)
- \int_{\sigma^l} f(z)d\rho(z).
\eeq
\end{enumerate}
\end{lemma}

\begin{proof}
 \begin{enumerate}
  \item 
Having in mind \eqref{intrep}, we will show as a first step that $\int_0^x \frac{Y^{1/2}(z)}{G(z,\tau)}d\tau$ is purely imaginary as $z\to E_{2j}$, with $z\in\sigma$, (the case $z\to E_{2j-1}$ can be handled in the same way).
For fixed $x\in\R$ we can separate the interval $[0,x]$ into smaller intervals $[0,x_1]\cup[x_1,x_2]\cup \dots \cup[x_k,x]$ such that $\mu_j(x_l)\in\{E_{2j-1},E_{2j}\}$ and $\mu_j(x_l)\not = \mu_j(x_{l+1})$. Assuming $\mu_j(x_l)=E_{2j-1}$ and $\mu_j(x_{l+1})=E_{2j}$, and setting 
\beq\nn
\tilde Y_j(z,x):=\sqrt{-(z-E_0)}\sqrt{z-E_{2j-1}}\prod_{j\not=l}\left(\frac{z-E_{2l-1}}{z-\mu_l(x)}\frac{z-E_{2l}}{z-\mu_l(x)}\right)^{1/2},
\eeq
where $\tilde Y_j(z,x)$ is bounded for any $z$ inside the j'th gap,
we can conclude 
\begin{align} \nn
& \int_{x_l}^{x_{l+1}}\frac{Y^{1/2}(z)}{G(z,\tau)}d\tau  =\sqrt{z-E_{2j}}\int_{x_l}^{x_{l+1}}\frac{\tilde Y_{j}(z,\tau)}{z-\mu_j(\tau)}d\tau \\ \nn
& = \sqrt{z-E_{2j}}\Big(\int_{x_l}^{x_{l+1}}\frac{\tilde Y_{j}(\mu_j(\tau),\tau)}{z-\mu_j(\tau)}d\tau +\int_{x_l}^{x_{l+1}}\frac{\tilde Y_{j}(z,\tau)-\tilde Y_{j}(\mu_j(\tau),\tau)}{z-\mu_j(\tau)}d\tau\Big) \\ \nn
& = \sqrt{z-E_{2j}}\Big( \int_{x_l}^{x_{l+1}} \frac{-\frac{d\mu_j(\tau)}{d\tau}}{2\sigma_j(\tau) \sqrt{\mu_j(\tau)-E_{2j}} (z-\mu_j(\tau))}d\tau+ \\ \nn 
& \qquad +\int_{x_l}^{x_{l+1}}\frac{d}{dz}\tilde Y_{j}(z,\tau) \vert_{ z=\zeta_j(\tau)}d\tau \Big),
\end{align}
where $\zeta_j(\tau)\in (\mu_j(\tau),z)$. Note that the function $\frac{d}{dz}\tilde Y_{j}(z,\tau)$ is uniformly bounded for $z\in[E_{2j-1}+\varepsilon,E_{2j}+\varepsilon]$ for some $\varepsilon>0$ and that $\frac{\tilde Y_{j}(z,\tau)-\tilde Y_{j}(\mu_j(\tau),\tau)}{z-\mu_j(\tau)}$ is uniformly bounded for $\mu_j\in[E_{2j-1}, E_{2j-1}+\varepsilon]$ and $z$ near $E_{2j}$, which yields
\beq
\sqrt{z-E_{2j}}\int_{x_l}^{x_{l+1}}\frac{d}{dz}\tilde Y_{j}(z,\tau) \vert_{ z=\zeta_j(\tau)}d\tau=O\left(\sqrt{z-E_{2j}}\right).
\eeq

On each of the intervals $[x_l,x_{l+1}]$ the function $\sigma_j(x)$ is constant and therefore  
\begin{align} 
&\sqrt{z-E_{2j}}\left( \int_{x_l}^{x_{l+1}} \frac{-\frac{d\mu_j(\tau)}{d\tau}}{2\sigma_j \sqrt{\mu_j(\tau)-E_{2j}} (z-\mu_j(\tau))}d\tau\right)\\ \nn
& = \sqrt{z-E_{2j}}\left(-\int_{\mu_j(x_l)}^{\mu_j(x_{l+1})} \frac{1}{2\sigma_j \sqrt{y-E_{2j}} (z-y)}dy \right) \\ \nn
& =\sqrt{E_{2j}-z}\left(\int_{\sqrt{E_{2j}-E_{2j-1}}}^{0} \frac{1}{\sigma_j (z-E_{2j}+s^2)}ds\right) \\ \nn
& = -\sigma_j\I\arctan\left( \frac{\sqrt{E_{2j}-E_{2j-1}}}{\sqrt{z-E_{2j}}}\right).
\end{align}

A close look shows that the same method can be applied to compute $\int_{0}^{x_1}\frac{Y^{1/2}(z)}{G(z,\tau)}d\tau$ and   $\int_{x_k}^{x}\frac{Y^{1/2}(z)}{G(z,\tau)}d\tau$ \\
This implies that $\int_{x_l}^{x_{l+1}}\frac{Y^{1/2}(z)}{G(z,\tau)}d\tau \to -\frac{1}{2}\sigma_j\I\pi$ as $z\to E_{2j}$ and thus  $\lim_{z\to E_{2j}}\int_0^x \frac{Y^{1/2}(z)}{G(z,\tau)}d\tau \in\I\R$. 

In more detail one obtains that 
\beq\label{intev}
\lim_{z\to E}\exp\left(\int_0^x \frac{Y^{1/2}(z)}{G(z,\tau)}d\tau\right)=\begin{cases}
 \pm 1, & \mu_j(0)\not =E, \mu_j(x)\not=E, \\
\pm 1, & \mu_j(0)=E, \mu_j(x)=E, \\
\pm\I, & \mu_j(0)=E, \mu_j(x)\not=E, \\
\pm\I, & \mu_j(0)\not=E, \mu_j(x)=E, \end{cases}
\eeq
for any $E\in\partial \sigma$.\\
For the product term in \eqref{intrep}, we have the following estimate
\begin{align}\nn
\frac{z-\mu_j(x)}{z-\mu_j(0)}& \exp\Big(-\frac{2}{\beta}\sum_{j=1}^\infty (E_{2j}-E_{2j-1})\Big)
\leq\prod_{l=1}^\infty \frac{z-\mu_l(x)}{z-\mu_l(0)}\\ \nn
&\leq \frac{z-\mu_j(x)}{z-\mu_j(0)}\exp\Big(\frac{2}{\beta}\sum_{j=1}^\infty (E_{2j}-{E_{2j-1}})\Big),
\end{align}
for $z$ inside the interval $[E_{2j-1}+\varepsilon,E_{2j}+\varepsilon]\cap\sigma$ for some $\varepsilon>0$, where $\beta =\min_{l,j;l\not = j}\{\vert E_{2j}-E_{2l}\vert, \vert E_{2j-1}-E_{2l}\vert \}$. This finishes the proof of \eqref{behpsi}, as \eqref{1.10} follows directly from \eqref{1.8}.

To show that a similar estimate holds for $\psi_\pm^\prime (z,x)$ as $z\to E\in\hat M$, consider, using \eqref{intrep},
\beq
\psi_\pm^\prime(z,x)=m_\pm(z,x)\psi_\pm(z,x).
\eeq
By the investigations of $\psi_\pm(z,x)$ from before, it suffices to analyze $m_\pm(z,x)$, given by \eqref{mfunc}, which has the following representation
\beq
m_\pm(z,x)=\sum_{n=1}^\infty \frac{\sigma_n(x) Y^{1/2}(\mu_n(x))}{\frac{d}{dz} G(\mu_n(x),x)(z-\mu_n(x))}\pm \frac{Y^{1/2}(z)}{G(z,x)}.
\eeq
For $j\not=n$ we have 
\beq\label{est1}
\left\vert \frac{Y^{1/2}(\mu_n(x))}{\frac{d}{dz}G(\mu_n(x),x)(E_{2j}-\mu_n(x))}\right\vert \leq C_1\frac{\sqrt{E_{2n}-E_0}}{\beta}(E_{2n}-E_{2n-1}),
\eeq
where $C_1:=\exp\Big(\frac{1}{\beta}\sum_{j=1}^\infty (E_{2j}-E_{2j-1})\Big)$.
Thus, using Hypothesis~\ref{Hyp}, we obtain that 
$\sum_{n\not =j} \frac{\sigma_j(x) Y^{1/2}(\mu_j(x))}{\frac{d}{dz} G(\mu_j(x),x)(z-\mu_j(x))}$ converges uniformly and is uniformly bounded with respect to $x$. \\
For $j=n$, we obtain 
\beq\label{est2}
\left\vert \frac{ Y^{1/2}(\mu_j(x))}{\frac{d}{dz} G(\mu_j(x),x)(z-\mu_j(x))}\right\vert 
\leq C_1 \frac{\sqrt{E_{2j}-E_0}\sqrt{(\mu_j(x)-E_{2j-1})(E_{2j}-\mu_j(x))}}{z-\mu_j(x)},
\eeq
and analogously
\beq\label{est3}
\left\vert \frac{Y^{1/2}(z)}{G(z,x)}\right\vert\leq C_1\frac{\sqrt{z-E_0}\sqrt{(z-E_{2j-1})(z-E_{2j})}}{z-\mu_j(x)}.
\eeq
Multiplying now \eqref{est1}, \eqref{est2}, and \eqref{est3} by $\psi_\pm(z,x)$, using the estimates obtained when proving \eqref{behpsi}, and letting $z\to E_{2j}$ with $z\in\sigma$ finishes the proof. 
\item
For a proof we refer to \cite[Sect. 9.5]{CL}, \cite[Thm 7.1.1]{L}, and \cite[Sect. 3.8]{T}.
 \end{enumerate}

\end{proof}

\section{Derivation of the integral equations for the transformation operators and estimates}

Consider the equation 
\beq
\left(-\frac{d^2}{dx^2}+q(x)\right)y(x)=zy(x), \quad z\in\C,
\eeq
with a potential $q(x)$ satisfying (\ref{moment}). Suppose that this equation has two solutions $\phi_\pm(z,x)$, which we will call the Jost solutions, which, for fixed $z$, are asymptotically close to the background Weyl solutions  $\psi_{\pm}(z,x)$ defined in (\ref{1.5}) as $x\to\pm\infty$. Set 
\begin{equation}\label{A.1}
J(z,x,y)=\frac{\psi_+(z,y)\psi_-(z,x) -
\psi_+(z,x)\psi_-(z,y)}{W(\psi_+(z),\psi_-(z))}
\end{equation}
and 
\begin{equation}\label{A.2}
\tilde q(x)=q(x) - p(x).
\end{equation}

It will be shown that the Jost solutions satisfy the following Volterra integral equations
\begin{equation}\label{A.3}
\phi_{\pm}(z,x) =\psi_{\pm}(z,x) -
\int_x^{\pm\infty}J(z,x,y)\tilde q(y)\phi_{\pm}(z,y)d y.
\end{equation}
We will also show that solutions of this form permit the following representation
\beq\label{asymprep}
\phi_{\pm}(z,x)=\psi_{\pm}(z,x)\pm\int_x^{\pm \infty} K_{\pm}(x,y)\psi_{\pm}(z,y)dy,
\eeq
where $K_\pm(x,y)$ are real-valued functions. 
For simplicity and because both representations can be obtained using the same techniques we will investigate the $+$ case only.

Assume that there exist $K_+(x,y)$ with $K_+(x,.)\in L^2(\R)$ and $K_+(x,y)=0$ for $ y< x$, such that $\phi_+(z,x)$ can be represented by (\ref{asymprep}).
Then substituting (\ref{asymprep}) into (\ref{A.3}), multiplying it with $\psi_-(z,x)g(z)$, integrating over the set $\sigma^{u,l}$, using the identity (\ref{1.14}), and taking into account that $K_+(x,y)=0$, $ x >  y$, we obtain 
\begin{equation}\label{A.4}
K_+(x,s) +\int_x^{\infty} dy\,
\tilde q(y)\oint_{\sigma}J(\lambda,x,y)\psi_+(\lambda,y)\psi_-(\lambda,s)
d\rho(\lambda)\qquad\qquad
\end{equation}
\[
\qquad +\int_x^{\infty}d y\,\tilde q(y)\int_y^{\infty} dt\,
K_+(y,t)\oint_{\sigma}J(\lambda,x,y)\psi_+(\lambda,t)\psi_-(\lambda,s)d\rho(\lambda)= 0.
\]

Set
\begin{equation}\label{A.5}
\Gamma_+(x,y,t,s)=\oint_{\sigma}
\psi_+(\lambda,x)\psi_-(\lambda,y)\psi_+(\lambda,t)\psi_-(\lambda,s) g(\lambda)d\rho(\lambda),
\end{equation}
where the integral has to be understood in the principal value sense
\footnote{
For informations about principal value integrals, we refer to \cite[Chap. 2]{Mu}.}.

Then substituting \eqref{1.62}, \eqref{1.12}, \eqref{A.1}, and \eqref{A.5} into \eqref{A.4} we obtain
\begin{align}\label{A.6}
&K_+(x,s) +\int_x^{\infty}\left(\Gamma_+(x,y,y,s) -
\Gamma_+(y,x,y,s)\right)\tilde q(y)\,d y\\
&+\int_x^{\infty}d y \,\tilde q(y)\int_y^{\infty}K_+(y,t)\left(\Gamma_+(x,y,t,s)
- \Gamma_+(y,x,t,s)\right)d t=0.\notag
\end{align}

A simple calculation using (\ref{1.8}) and (\ref{1.10}) shows that (\ref{A.5}) satisfies
\begin{equation}\label{A.7}
\overline{\Gamma_+(x,y,t,s)}=-\Gamma_+(y,x,s,t).
\end{equation}
Combining (\ref{intrep}), (\ref{1.88}), and (\ref{1.12}), one obtains that the only poles of the integrand are given at the band edges. Using a series of contour integrals and applying Jordan's lemma (see for example \cite{STSH}), it turns out to be necessary to investigate the following series 
\beq\label{3.2}
D_+(x,y,r,s)=-\frac{1}{4}\sum_{E\in\partial \sigma} f_+(E,x,y,r,s),
\eeq
where, for any $E\in\partial\sigma$, we denote 
\beq\label{3.1}
f_+(E,x,y,r,s)=\lim_{z\to E}
\frac{G(z,0)^2}{\frac{d}{dz}Y(z)}
\psi_+(z,x)\psi_-(z,y) \psi_+(z,r) \psi_-(z,s).
\eeq

\begin{lemma}\label{D}
The series $D_+(x,y,r,s)$ defined by (\ref{3.2}) and (\ref{3.1}) converges, is  continuous, and uniformly bounded with respect to all variables $(x,y,r,s)\in\R^4$.
\end{lemma}

\begin{proof}
First note that by (\ref{intrep}) we have 
\begin{align}\nn
f_+(E,x,y,r,s)& =\frac{(G(E,x) G(E,y) G(E,r) G(E,s))^{1/2}}{\frac{d}{dz}Y(E)}  \\ \nn
& \quad \times \lim_{z\to E} \exp\left(  \int_y^x \frac{Y^{1/2}(z)}{G(z,\tau)}d\tau + \int_s^r \frac{Y^{1/2}(z)}{G(z,\tau)}d\tau\right),
\end{align}
for some $E\in\partial\sigma$,
where the limit is taken from inside the spectrum.
For computing the integral terms we refer to the proof of Lemma~\ref{lem1.1}.
We will now investigate 
\begin{align}\nn
M_+(E_{2l},x,y,r,s)& = \frac{(G(E_{2l},x) G(E_{2l},y) G(E_{2l},r) G(E_{2l},s))^{1/2}}{\frac{d}{dz}Y(E_{2l})}\\
& =-\frac{((E_{2l}-\mu_l(x))(E_{2l}-\mu_l(y))(E_{2l}-\mu_l(r))(E_{2l}-\mu_l(s)))^{1/2}}{(E_{2l}-E_0)(E_{2l}-E_{2l-1})} \\ \nn
& \quad\times \prod_{j=1, j\not = l}^\infty \frac{((E_{2l}-\mu_j(x))(E_{2l}-\mu_j(y))(E_{2l}-\mu_j(r))(E_{2l}-\mu_j(s)))^{1/2}}{(E_{2l}-E_{2j-1})(E_{2l}-E_{2j})},
\end{align}
which can be estimated as follows
\begin{align}\nn
 \prod_{j=1,j\not = l}^\infty & \vert \frac{(E_{2l}-\mu_j(x))}{(E_{2l}-E_{2j})}  \frac{(E_{2l}-\mu_j(y))}{(E_{2l}-E_{2j-1})}\vert 
\leq \prod_{j=1}^{l-1} \frac{(E_{2l}-\mu_j(x))}{(E_{2l}-E_{2j})}\prod_{j=l+1}^{\infty}\frac{(E_{2l}-\mu_j(y))}{(E_{2l}-E_{2j-1})}\\ \nn
& \leq \exp\left(\sum_{j=1}^{l-1} \log\Big(1+\frac{E_{2j}-\mu_j(x)}{E_{2l}-E_{2j}}\Big) +\sum_{j=l+1}^\infty \log\Big(1+\frac{\mu_j(y)-E_{2j-1}}{E_{2j}-E_{2l}}\Big)\right) \\ \nn
& \leq \exp \Big(\frac{1}{\beta}\sum_{j=1}^\infty (E_{2j}-E_{2j-1}) \Big)<\infty,
\end{align}
where we used $\log(1+x)\leq x$ for $x>0$ and $\beta =\min_{l,j;l\not = j}\{\vert E_{2j}-E_{2l}\vert, \vert E_{2j-1}-E_{2l}\vert \}$.
Moreover, 
\begin{align} \label{conv} 
\frac{((E_{2l}-\mu_l(x))(E_{2l}-\mu_l(y))(E_{2l}-\mu_l(r))(E_{2l}-\mu_l(s)))^{1/2}}{(E_{2l}-E_0)(E_{2l}-E_{2l-1})} \leq \frac{(E_{2l}-E_{2l-1})}{(E_{2l}-E_0)}.
\end{align}
This implies
\beq\label{estM}
0\leq\vert f_+(E_k,x,y,r,s)\vert = \vert M_+(E_{k},x,y,r,s)\vert \leq C_1\frac{E_{2l}-E_{2l-1}}{E_{2l}-E_0},
\eeq
where $k\in \{2l-1,2l\}$ and $C_1:=\exp\Big(\frac{1}{\beta}\sum_{j=1}^\infty (E_{2j}-E_{2j-1})\Big)$, and therefore our series 
converges and hence $D_+(x,y,r,s)$ is well defined and uniformly bounded with respect to all variables.

The continuity in all variables follows immediately by using 
\eqref{intrep} and \eqref{intev}. 
\end{proof}

\begin{lemma}
The function $D_+(x,y,r,s)$, defined through \eqref{3.2} and \eqref{3.1}, has first partial derivatives, which are uniformly bounded on $\R^4$.
\end{lemma}

\begin{proof}
Consider the integral representation (\ref{intrep}) of the background Weyl solutions, 
then the derivative with respect to $x$ is given by 
\beq
\psi_+^\prime(z,x)=\left( \frac{H(z,x)+ Y^{1/2}(z)}{G(z,x)} \right) \left( \frac{G(z,x)}{G(z,0)}\right)^{1/2}\exp\left( \int_0^x \frac{Y^{1/2}(z)}{G(z,x)} \right).
\eeq
Note that $G(z,0)^{1/2}\frac{d}{dx}\psi_+(z,x)$, by Lemma~\ref{lem1.1} has neither poles nor square root singularities at the band edges $E\in\partial\sigma$, which allows us to pass to the limit in the following expression
\begin{align}
\lim_{z\to E}(z-E)\frac{H(z,x) + Y^{1/2}(z)}{G^{1/2}(z,x)} \frac{(G(z,y)G(z,r)G(z,s))^{1/2}}{Y(z)}.
\end{align}
Therefore, we will slightly abuse notation by omitting the limit and replacing $z$ by $E$. 
W.l.o.g. we will assume that $E=E_{2n}$ (the case $E=E_{2n-1}$ can be treated similarly).
Using \eqref{repHG}, we have 
\begin{align}\nn
\frac{H(E_{2n},x)}{G^{1/2}(E_{2n},x) } & \frac{(G(E_{2n},y)G(E_{2n},r)G(E_{2n},s))^{1/2}}{\frac{d}{dz}Y(E_{2n})}\\ \nn
&  = \sum_{j=1}^\infty \frac{\sigma_j(x) Y^{1/2}(\mu_j(x))}{\frac{d}{dz} G(\mu_j(x),x)(E_{2n}-\mu_j(x))}M_+(E,x,y,r,s).
\end{align}
Due to \eqref{estM} we will at first consider $\sum_{j\not =n}\frac{Y^{1/2}(\mu_j(x))}{\frac{d}{dz} G(\mu_j(x),x)(E_{2n}-\mu_j(x))}$, which can be done using the same ideas as in Lemma~\ref{D}.
Namely, for $j\not =n$,
\beq
0 \leq \left\vert \frac{Y^{1/2}(\mu_j(x))}{\frac{d}{dz} G(\mu_j(x),x)(E_{2n}-\mu_j(x))} \right\vert 
\leq \frac{1}{\beta} \sqrt{(E_{2j}-E_0)}(E_{2j}-E_{2j-1})C_1^{1/2},
\eeq
which implies that the corresponding sum converges. \\
For $j=n$, there are two cases to distinguish:
\begin{enumerate}
\item
If $\mu_n(x)=E_{2n}$, then 
\begin{align}\nn
 0& \leq \left\vert \frac{Y^{1/2}(E_{2n})}{G(E_{2n},x)}M_+(E_{2n},x,y,r,s)\right\vert \\ \nn
&\leq C_1^{3/2}\frac{((E_{2n}-\mu_n(y))(E_{2n}-\mu_n(r))(E_{2n}-\mu_n(s)))^{1/2}}{\sqrt{(E_{2n}-E_0)(E_{2n}-E_{2n-1})}}\\ \nn
&\leq C_1^{3/2} \frac{E_{2n}-E_{2n-1}}{\sqrt{E_{2n}-E_0}}.
\end{align}

\item
If $\mu_n(x)\not = E_{2n}$, we have 
\begin{align}\nn
0 & \leq \left\vert \frac{Y^{1/2}(\mu_n(x))}{\frac{d}{dz} G(\mu_n(x),x)(E_{2n}-\mu_n(x))}M_+(E_{2n},x,y,r,s)\right\vert\\ \nn
& \leq C_1^{3/2}(\mu_n(x)-E_0)^{1/2}\frac{((\mu_n(x)-E_{2n-1})(E_{2n}-\mu_n(y)))^{1/2}}{(E_{2n}-E_0)}\\ \nn
& \qquad\times\frac{((E_{2n}-\mu_n(r))(E_{2n}-\mu_n(s)))^{1/2}}{(E_{2n}-E_{2n-1})}\\ \nn
& \leq C_1^{3/2} \frac{(E_{2n}-E_{2n-1})}{\sqrt{E_{2n}-E_0}}.
\end{align}
\end{enumerate}
Next we consider 
\beq
\frac{Y^{1/2}(E_{2n})}{G^{1/2}(E_{2n},x)} \frac{(G(E_{2n},y)G(E_{2n},r)G(E_{2n},s))^{1/2}}{\frac{d}{dz}Y(E_{2n})},
\eeq
which can be investigated as before. 
\begin{enumerate}
\item
If $\mu_n(x)=E_{2n}$, then 
\begin{align}\nn
0& \leq \left\vert  \frac{Y^{1/2}(E_{2n})}{G(E_{2n},x)}M_+(E_{2n},x,y,r,s)\right\vert\\ \nn  
&\leq C_1^{3/2}\frac{ ((E_{2n}-\mu_n(y))(E_{2n}-\mu_n(r))(E_{2n}-\mu_n(s)))^{1/2}}{\sqrt{(E_{2n}-E_0)(E_{2n}-E_{2n-1})}}\\ \nn
 & \leq C_1^{3/2}\frac{(E_{2n}-E_{2n-1})}{\sqrt{(E_{2n}-E_0)}}.
\end{align}
\item
If $\mu_n(x)\not = E_{2n}$, 
\beq\nn 
 \frac{Y^{1/2}(E_{2n})}{G(E_{2n},x)}M_+(E_{2n},x,y,r,s)=0,
\eeq 
because 
\beq
\lim_{z\to E_{2n}}\frac{(-(z-E_{2n-1})(z-E_{2n})(z-\mu_n(y))(z-\mu_n(r))(z-\mu_n(s)))^{1/2}}{\sqrt{(z-E_0)(z-\mu_n(x))}(z-E_{2n-1})} =0.
\eeq
\end{enumerate}
This finishes the proof, as all our estimates are uniformly with respect to $x,y,r,s\in\R$ and 
\begin{equation}
 0\leq \vert f_{+,x}(E_{k},x,y,r,s)\vert \leq C\frac{(E_{2n}-E_{2n-1})}{\sqrt{(E_{2n}-E_0)}},
\end{equation}
where $k\in\{2n-1,2n\}$ and $C$ denotes a constant independent of $n$.
\end{proof}

\begin{lemma}\label{lemA.1}
The kernels $K_\pm(x,s)$ of the transformation operators satisfy
the integral equation
\begin{align}\nn
K_\pm(x,s)
&= - 2\int_{\frac{x+s}{2}}^{\pm\infty}\tilde q(y)D_\pm(x,y,y,s)d y\\ \label{inteq:K}
&\quad
 \mp 2\int_x^{\pm \infty}d y\int_{s+x-y}^{s+y-x}
D_\pm(x,y,r,s)K_\pm(y,r)\tilde q(y)\,d r,\quad\pm s >\pm x,
\end{align}
where $D_\pm(x,y,r,s)$ are defined by (\ref{3.2}).
In particular, 
\beq\label{Kxx}
K_\pm(x,x)=\pm\frac{1}{2}\int_x^{\pm\infty} (q(s)-p(s)) ds.
\eeq
\end{lemma}

\begin{proof}
Suppose that $(x-y+r-s)>0$, where $x,y,r,s$ are considered fixed parameters, and take a series of closed contours $C_n$ consisting of a circular arc $R_n$ centered at the origin with radius $(E_{2n}+E_{2n+1})/2$, together with some parts wrapping around each of the first $n+1$ bands of the spectrum $\sigma$, but not intersecting it, as indicated in figure~\ref{C_n}. Then we can consider the following series of contour integrals
\begin{equation}
 I_n(x,y,r,s):=\frac{1}{2\pi\I}\int_{C_n} \psi_+(z,x)\psi_-(z,y)\psi_+(z,r)\psi_-(z,s)g(z)^2dz.
\end{equation}

\begin{figure}
\begin{picture}(7,5.2)
\put(1,2.6){\line(1,0){0.5}}
\put(2.3,2.6){\line(1,0){1.3}}
\put(4.2,2.6){\line(1,0){1.8}}

\put(4.5,2.75){\vector(1,0){0.1}}
\put(4.5,2.45){\vector(-1,0){0.1}}
\put(3.0,2.75){\vector(1,0){0.1}}
\put(3.0,2.45){\vector(-1,0){0.1}}
\put(1.2,2.75){\vector(1,0){0.1}}
\put(1.2,2.45){\vector(-1,0){0.1}}

\curve(4.4,0.2,  4.6,0.8,  4.86,2.0,   4.892,2.3, 4.8,2.45,  4.2,2.45,   4.05,2.6,   4.2,2.75, 4.8, 2.75,  4.892,2.9,   4.86,3.2,  4.6,4.4,  4.4,5.)
\curve( 2.15,2.6, 2.3,2.75,  2.6,2.75, 3.6,2.75, 3.75,2.6, 3.6,2.45, 2.6,2.45, 2.3,2.45, 2.15,2.6)
\curve(1.,2.75, 1.2,2.75, 1.5,2.75, 1.65,2.6, 1.5,2.45, 1.2,2.45, 1.,2.45) 

\put(4.2,2.2){$\scriptstyle E_{2n}$}
\put(3.0,2.2){$\scriptstyle E_{2n-1}$}
\put(2.2,2.2){$\scriptstyle E_{2n-2}$}
\put(1.3,2.2){$\scriptstyle E_{2n-3}$}
\end{picture}
\caption{Contour $C_n$}
\label{C_n}
\end{figure}
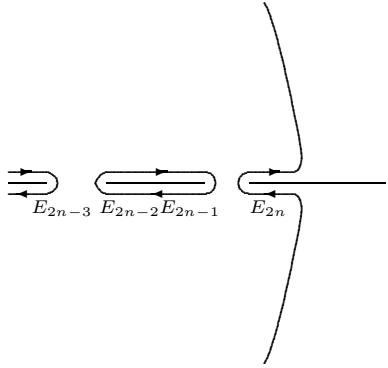

Moreover, consider $L_n$ the integral along the loop around the n'th band of the spectrum and set 
\begin{equation}
 S_n= \frac{1}{2\pi\I}\oint_{[E_{2n-2},E_{2n-1}]} \psi_+(z,x)\psi_-(z,y)\psi_+(z,r)\psi_-(z,s)g(z)^2dz,
\end{equation}
then we have 
\begin{equation}
 L_n=L_n-S_n+S_n=S_n, 
\end{equation}
where we used that $L_n-S_n=0$, by the continuity of the integrand up to the real axis and Cauchy's theorem.

Thus, shrinking the loops around each band of the spectrum in $C_n$ to an integral along each band of the spectrum in $C_n$, and denoting by $\tilde I_n(x,y,r,s)$ the resulting series of contour integrals, we obtain 
\begin{align}
\tilde I_n(x,y,r,s)&=\frac{1}{2\pi\I}\int_{R_n} \psi_+(z,x)\psi_-(z,y)\psi_+(z,r)\psi_-(z,s)g(z)^2dz\\ \nn 
&\quad - \frac{1}{4}\sum_{l=0}^{2n}f_+(E_l,x,y,r,s),
\end{align}
where we used Sokhotsky--Plemelj's formula (see e.g. \cite[Chap. 2]{Mu}).

On the circle $R_n$ we have the following asymptotic behavior as $n\to\infty$ and therefore $z\to\infty$,
\beq
g^{1/2}(z)\psi_+(z,x)=\E^{\I\sqrt{z}x(1+O(\frac{1}{z}))}O\left(\frac{1}{z^{1/4}}\right),
\eeq
where we used Lemma~\ref{lempsi} and the fact that these asymptotics are valid as long as we are outside a small neighborhood of the gaps. This yields 
\beq
\psi_+(z,x)\psi_-(z,y)\psi_+(z,r)\psi_-(z,s)g(z)^2= \E^{\I\sqrt{z}(x-y+r-s)(1+O(\frac{1}{z}))}O\left( \frac{1}{z} \right),
\eeq
as $z\to\infty$. 

Hence one can apply Jordan's lemma to conclude that the contribution of the circle $R_n$ vanishes as $n\to\infty$. 

Thus letting $n\to\infty$ the series of integrals $\tilde I_n(x,y,r,s)$ converges to 
\begin{equation}\label{A.10}
\Gamma_+(x,y,r,s)=D_+(x,y,r,s),\quad\mbox{for}\quad
(x-y+r-s)>0.
\end{equation}

Note that $f_+(E,x,y,r,s)$ is real for any $E\in\partial\sigma$, because \eqref{intev} and $G(E,x)=0$, if $\mu_j(x)=E$, imply that $f_+(E,x,y,r,s)=0$.
Moreover $f_+(E,x,y,r,s)=f_+(E,y,x,s,r)$. Thus $D_+(x,y,r,s)$ is also real and satisfies 
\beq\label{A.11}
D_+(x,y,r,s)=D_+(y,x,s,r).
\eeq

Now let $(x-y+r-s)<0$, that is $-(x-y+r-s)>0$. Then \eqref{A.7}, \eqref{A.10}, and \eqref{A.11} imply 
\beq
\Gamma_+(x,y,t,s)=-\overline{\Gamma_+(y,x,s,r)}=- \overline{D_+(x,y,r,s)}=-D_+(x,y,r,s).
\eeq

Therefore,
\beq
\Gamma_+(x,y,r,s)=D_+(x,y,r,s)\sign(x-y+r-s).
\eeq

Combining all the facts, the domain, where the first integrand in (\ref{A.6}) does not vanish, is given by 
\begin{equation}\label{A.13}
\sign(x-s)=-\sign(2y-x-s),\quad  s> x.
\end{equation}
In the second integral the domain of integration is
\begin{equation}\label{A.121}
\sign(x-y+t-s)=-\sign(y-x+t-s),
\text{ with}\quad s> x,\quad t> y> x.
\end{equation}

Solving (\ref{A.13}) and (\ref{A.121}), proves \eqref{inteq:K}.

Setting now $s=x$ in (\ref{inteq:K}), the second summand vanishes, because we set $K_+(y,r)=0$ for $ r<  y$. Hence 
\beq\label{A.16}
K_+(x,x)=-2\int_x^{\infty} \tilde q(y)D_+(x,y,y,x)dy.
\eeq

Thus we obtain 
\begin{align}
D_+(x,y,y,x)
& =\frac{1}{4}\sum_{E\in\partial\sigma} \Res_E \frac{1}{z-E_0}\prod_{j=1}^\infty \frac{(z-\mu_j(x))(z-\mu_j(y))}{(z-E_{2j-1})(z-E_{2j})}\\ \nn
& =\frac{1}{4}\sum_{l=0}^\infty \lim_{z\to E_l}\frac{z-E_l}{z-E_0}\prod_{j=1}^\infty \frac{(z-\mu_j(x))(z-\mu_j(y))}{(z-E_{2j-1})(z-E_{2j})},
\end{align}
and we already know that this function is bounded by Lemma~\ref{D}. 
Considering now the following sequence 
\begin{align}
D_{+,n}(x,y,y,x)
& =\frac{1}{4}\sum_{E\in\partial\sigma} \Res_E \frac{1}{z-E_0}\prod_{j=1}^n \frac{(z-\mu_j(x))(z-\mu_j(y))}{(z-E_{2j-1})(z-E_{2j})}\\ \nn
& =\frac{1}{4}\sum_{l=0}^{2n}\lim_{z\to E_l} \frac{z-E_l}{z-E_0}\prod_{j=1}^n \frac{(z-\mu_j(x))(z-\mu_j(y))}{(z-E_{2j-1})(z-E_{2j})},
\end{align}
which corresponds to the case where we only have $n$ gaps and crossed out all the other ones. 
We now estimate  
\begin{align}\label{estDN}
 \vert D_+(x,y,y,x)-D_{n,+}(x,y,y,x)\vert 
& \leq \frac{1}{4} \sum_{l=0}^{2n} \lim_{z\to E_l} \left\vert \frac{z-E_l}{z-E_0}\prod_{j=1}^n \left(\frac{z-\mu_j(x)}{z-E_{2j-1}}\frac{z-\mu_j(y)}{z-E_{2j}}\right)\right\vert \\ \nn 
&\qquad \times \left\vert \left(\prod_{j=n+1}^\infty \left(\frac{z-\mu_j(x)}{z-E_{2j-1}}\frac{z-\mu_j(y)}{z-E_{2j}} \right)-1 \right)\right\vert \\ \nn
& + \frac{1}{4} \sum_{l=2n+1}^\infty \lim_{z\to E_l} \left\vert \frac{z-E_l}{z-E_0}\prod_{j=1}^\infty \left( \frac{z-\mu_j(x)}{z-E_{2j-1}}\frac{z-\mu_j(y)}{z-E_{2j}}\right)\right\vert,
\end{align}
using the same techniques as in the proof of Lemma~\ref{D}.
We fix $z=E_{2l}$ (the case $z=E_{2l-1}$ can be handled analogously). If $l<n$, we have 
\begin{align}\nn
 \lim_{z\to E_{2l}}& \left\vert \frac{z-E_{2l}}{z-E_0}\prod_{j=1}^n\left(\frac{z-\mu_j(x)}{z-E_{2j-1}}\frac{z-\mu_j(y)}{z-E_{2j}}\right)\right\vert\\ \nn 
 & =\left\vert\frac{E_{2l}-\mu_l(x)}{E_{2l}-E_0}\frac{E_{2l}-\mu_l(y)}{E_{2l}-E_{2l-1}}\prod_{j=1, j\not = l}^n\frac{E_{2l}-\mu_j(x)}{E_{2l}-E_{2j-1}}\frac{E_{2l}-\mu_j(y)}{E_{2l}-E_{2j}}\right\vert \\ \nn
& \leq\frac{E_{2l}-E_{2l-1}}{E_{2l}-E_0}\prod_{j=1}^{l-1}\frac{E_{2l}-\mu_j(x)}{E_{2l}-E_{2j}}\prod_{j=l+1}^n\frac{E_{2l}-\mu_j(x)}{E_{2l}-E_{2j-1}} \\ \nn
& \leq \frac{E_{2l}-E_{2l-1}}{E_{2l}-E_0}\exp\Big(\frac{1}{\beta}\sum_{j=1}^n (E_{2j}-E_{2j-1})\Big)
\end{align}
and 
\begin{align}
\exp\Big(-\frac{1}{\beta}& \sum_{j=n+1}^\infty (E_{2j}-E_{2j-1})\Big)
 \leq\prod_{j=n+1}^\infty\frac{\mu_j(x)-E_{2l}}{E_{2j}-E_{2l}}\\ \nn
& \leq\prod_{j=n+1}^\infty\frac{\mu_j(x)-E_{2l}}{E_{2j-1}-E_{2l}}\frac{\mu_j(y)-E_{2l}}{E_{2j}-E_{2l}}\\ \nn
& \leq \prod_{j=n+1}^\infty \frac{\mu_j(x)-E_{2l}}{E_{2j-1}-E_{2l}}
\leq \exp\Big(\frac{1}{\beta}\sum_{j=n+1}^\infty (E_{2j}-E_{2j-1})\Big),
\end{align}
where the last estimate implies that the first sequence in \eqref{estDN} converges uniformly to zero as $n$ tends to $\infty$ as in the investigation of (\ref{quasiG}). 
Analogously, one can estimate the terms of the second sequence to obtain
\begin{equation}\nn
\lim_{z\to E_{2l}} \left\vert \frac{z-E_l}{z-E_0}\prod_{j=1}^\infty \left( \frac{z-\mu_j(x)}{z-E_{2j-1}}\frac{z-\mu_j(y)}{z-E_{2j}}\right)\right\vert
\leq \frac{E_{2l}-E_{2l-1}}{E_{2l}-E_{0}}\exp\Big(\frac{1}{\beta}\sum_{j=1}^\infty (E_{2j}-E_{2j-1})\Big).
\end{equation}
Thus also the second sequence in \eqref{estDN} converges uniformly to zero, as we are working in the Levitan class and hence $D_{n,+}(x,y,y,x)$ converges uniformly against $D_+(x,y,y,x)$.

 Moreover, it is known (see e.g. \cite{BET}), that 
\beq
D_{n,+}(x,y,y,x)=-\frac{1}{4}
\eeq
for each fixed $n$,
and hence we finally obtain
\beq
D_+(x,y,y,x)=\lim_{n\to\infty} D_{n,+}(x,y,y,x)=\lim_{n\to\infty}-\frac{1}{4}= -\frac{1}{4}.
\eeq

Therefore we can now conclude, using (\ref{A.16}), that 
\beq\label{Kqp}
K_\pm(x,x)=\pm\frac{1}{2}\int_x^{\pm\infty} (q(s)-p(s)) ds.
\eeq
\end{proof}

\begin{lemma}\label{lem:est}
Suppose (\ref{moment}), then (\ref{inteq:K}) has a unique solution $K_\pm(x,y)$ such that $K_\pm(x,y)$ has first-order partial derivatives with respect to both variables.
Moreover, for $\pm y \geq \pm  x$ the following estimates are valid
\begin{align}\label{A.21}
\left|K_\pm(x,y)\right|
&\leq C_\pm(x) Q_\pm(x+y),\\
\label{A.22}
\left|\frac{\partial K_\pm(x,y)}{\partial x}\right| +
\left|\frac{\partial K_\pm(x,y)}{\partial y}\right|
&\leq
C_\pm(x)\,\left(\left|\tilde q\left(\frac{x+y}{2}\right)\right|
+ Q_\pm\left(x+y\right)\right),
\end{align}
where
\beq
Q_\pm(x)=\pm\int_{\frac{x}{2}}^{\pm\infty} \vert \tilde q(s) \vert ds, \quad \tilde q(x)=q(x)-p(x),
\eeq
and $C_\pm(x)$ are positive continuous functions for $x\in\R$, which decrease as $x\to\pm\infty$ and depend on the corresponding background data and on $\pm\int_{2x}^{\pm\infty} Q_\pm(s)ds$.
\end{lemma}

\begin{proof}
Using the method of successive approximation one can prove existence and uniqueness of the solution $K_\pm(x,y)$ of (\ref{inteq:K}).
We restrict our considerations to the $+$ case. After the following change of variables 
\begin{equation}
 2\alpha:=s+r, \text{ }2\beta:=r-s,\text{ } 2u:=x+y, \text{ }2v:=y-x,
\end{equation}
\eqref{inteq:K} becomes
\begin{align}\nn
 H(u,v)  & = -2\int_u^\infty \tilde q(s)D_1(u,v,s)ds \\ 
& -4\int_u^{\infty} d\alpha \int_0^v \tilde q(\alpha-\beta) D_2(u,v,\alpha,\beta)H(\alpha,\beta) d\beta,
\end{align}
with
\begin{align}\nn
 H(u,v)=K_+(u-v,u+v), & \quad D_1(u,v,s)=D_+(u-v,s,s,u+v),\\
& D_2(u,v,\alpha,\beta)=D_+(u-v, \alpha-\beta, \alpha+\beta,u+v).
\end{align}
As the functions $D_1$ and $D_2$ are uniformly bounded with respect to all their variables by a constant $C$, we can apply the method of successive approximation to estimate $H(u,v)$, which yields 
\begin{equation}\label{estH}
 \vert H(u,v)\vert \leq C(u-v)Q_+(2u),
\end{equation}
where 
\beq\label{C}
C(u-v)=2C\exp\left(4C\int_{2u-2v}^\infty Q_+(x)dx\right).
\eeq

To obtain the second estimate \eqref{A.22} we recall that the partial derivatives with respect to all variables exist for $D_1$ and $D_2$ and that they are also bounded with respect to all variables.
Thus, using 
\begin{align}
& \frac{\partial H(u,v)}{\partial u}-2\tilde q(u)D_1(u,v,u)=\\ \nn
& =+4\int_0^v \tilde q(u-\beta)D_2(u,v,u,\beta)H(u,\beta)d\beta-2\int_u^\infty \tilde q(s)\frac{\partial D_1(u,v,s)}{\partial u}ds \\ \nn
& -4\int_u^\infty d\alpha\int_0^v\tilde q(\alpha-\beta)\frac{\partial D_2(u,v,\alpha,\beta)}{\partial u}H(\alpha,\beta)d\beta, 
\end{align}
and
\begin{align}
& \frac{\partial H(u,v)}{\partial v}= \\ \nn
& =-2\Big(\int_u^\infty \tilde q(s)\frac{\partial D_1(u,v,s)}{\partial v }ds+2\int_u^\infty \tilde q(\alpha-\beta)D_2(u,v,\alpha,v)H(\alpha,v)d\alpha \\ \nn
& +2\int_u^\infty d\alpha \int_0^v \tilde q(\alpha-\beta)\frac{\partial D_2(u,v,\alpha,\beta)}{\partial v }H(\alpha,\beta)d\beta\Big),
\end{align}
one obtains
\begin{align}
 & \vert \frac{\partial}{\partial u} H(u,v)\vert \leq C_1(u-v)(\vert \tilde q(u)\vert +Q_+(2u)),\\ \nn
& \vert \frac{\partial}{\partial v} H(u,v)\vert \leq C_1(u-v)(\vert \tilde q(u)\vert +Q_+(2u)),
\end{align}
where $C_1(u-v)$ is of the same type as \eqref{C} with a different constant $C_1$ depending on the background data. Using 
\begin{align}
\left\vert  \frac{dK_+(x,y)}{dx} \right\vert+\left\vert \frac{dK_+(x,y)}{dy} \right\vert \leq 
\left\vert\frac{dH(u,v)}{du} \right\vert+\left\vert \frac{dH(u,v)}{dv} \right\vert,
\end{align}
completes the proof.
\end{proof}

To finish the proof of Theorem~\ref{mainth}, we have to show the following:

\begin{lemma}
The functions 
\begin{equation}\label{eq:tiphi}
 \phi_\pm(z,x)= \psi_\pm(z,x)\pm\int_x^{\pm\infty}K_\pm(x,s)\psi_\pm(z,s)ds, 
\end{equation}
where $K_\pm(x,s)$ are defined by \eqref{inteq:K}, satisfy 
\begin{equation}\label{schrtiphi}
 \Big(-\frac{d^2}{dx^2}+q(x)\Big)\phi_\pm(z,x)=z\phi(z,x).
\end{equation}
\end{lemma}

\begin{proof}
 Again we will only consider  the $+$ case as the other one can be treated similarly and drop the $+$ whenever possible. 
On the one hand we obtain, using \eqref{Kxx}, that 
\begin{align}\nn
 \Big(-\frac{d^2}{dx^2}+q(x)\Big)\phi(z,x)& = -\psi^{\prime\prime}(z,x)+p(x)\psi(z,x)\\ \nn 
& +\frac{1}{2}\tilde q(x) \psi(z,x)+K(x,x)\psi^\prime(z,x)+K_x(x,x)\psi(z,x)\\ 
& +\int_x^\infty (q(x)K(x,s)-K_{xx}(x,s))\psi(z,s)ds,
\end{align}
and on the other hand, using that $\psi_\pm(z,x)$ are the background Weyl solutions, we have
\begin{align}\nn
 z\phi(z,x)& = z\psi(z,x)+K(x,x)\psi^\prime(z,x)-K_y(x,x)\psi(z,x)\\
&  +\int_x^\infty (p(s)K(x,s)-K_{ss}(x,s)) \psi(z,s)ds. 
\end{align}
Applying \eqref{Kxx} once more, we see that \eqref{schrtiphi} is satisfied if and only if 
\begin{equation}\label{diff K}
 \int_x^\infty (K_{xx}(x,s)-K_{ss}(x,s))\psi(z,s)ds=\int_x^\infty (q(x)-p(s))K(x,s)\psi(z,s)ds. 
\end{equation}
For proving this identity we use the integral equation \eqref{A.4} instead of \eqref{inteq:K} for $K(x,s)$, which yields for $x<s$
\begin{align}\nn
 K_{ss}(x,s)& = -\int_x^\infty dy \tilde q(y) \oint_\sigma J(\la,x,y)\psi_+(\la,y)\psi_-^{\prime\prime}(\la,s)d\rho (\la)\\ \nn 
& -\int_x^\infty dy \tilde q(y)\int_y^\infty dt K(y,t) \oint_\sigma J(\la,x,y) \psi_+(\la,t) \psi_-^{\prime\prime}(\la,s) d\rho(\la)\\ \nn 
&= (p(s)-p(x))K(x,s)-\int_x^\infty dy \tilde q(y) \oint_\sigma J_{xx}(\la,x,y)\psi_+(\la,y)\psi_-(\la,s)d\rho(\la)\\ \nn
& \quad -\int_x^\infty dy \tilde q(y) \int_y^\infty dt K(y,t) \oint_\sigma J_{xx}(\la,x,y) \psi_+(\la,t)\psi_-(\la,s)d\rho(\la)\\ \nn 
& = (p(s)-p(x))K(x,s) +K_{x,x}(x,s)-\tilde q(x) \oint_\sigma J_x(\la,x,x)\psi_+(\la,x)\psi_-(\la,s) d\rho(\la)\\  \nn 
& \quad -\tilde q(x) \int_x^\infty dt K(x,t)\oint_\sigma J_x(\la,x,x)\psi_+(\la,t)\psi_-(\la,s) d\rho(\la)\\ 
& = (p(s)-q(x))K(x,s) +K_{xx}(x,s).   
\end{align}
Here it should be noticed that $\oint_\sigma J(\la,x,y) \psi_+(\la,t)\psi_-^{\prime\prime}(\la,s) d\rho(\la)$ exists, because we can again estimate the sum of the absolute values of the residues by using that $-\psi_{\pm}^{\prime\prime}(z,x)+p(x)\psi_\pm(z,x)=z\psi_\pm(z,x)$ and the same techniques as in Lemma~\ref{D}. Analogously for $\oint_\sigma J_{xx}(\la,x,y)\psi_+(\la,t)\psi_-(\la,s)d\rho(\la)$. Thus \eqref{diff K} is fulfilled and therefore also \eqref{schrtiphi}. 

\end{proof}

\begin{lemma}
One has $K_\pm(x,y)\in L^2(\R_\pm, dy)$.
\end{lemma}
\begin{proof}
Using (\ref{moment}), we conclude
\begin{align}
\int_{\R}\vert K_\pm(x,y)\vert^2 dy& \leq \pm C_\pm(x)^2 \int_{x}^{\pm\infty} Q_\pm(x+y)^2 dy \\ \nn
&= C_\pm(x)^2Q_\pm(2x)\int_x^{\pm\infty}\int_{\frac{x+y}{2}}^{\pm\infty}\vert\tilde q(s)\vert dsdy\\ \nn
&=C_\pm(x)^2Q_\pm(2x)\int_x^{\pm\infty} (2s-2x)\vert \tilde q(s) \vert ds<\infty.
\end{align}
\end{proof}

It should also be noticed that for any function $f_\pm(x)\in L^2(\R_\pm)$, 
\beq
h_\pm(x)=f_\pm(x)\pm\int_x^{\pm\infty}K_\pm(x,y)f_\pm(y)dy\in L^2(\R_\pm).
\eeq
Thus as a consequence we obtain

\begin{corollary}
 The normalized Jost solutions $\phi_\pm(z,x)$ coincide with the Weyl solutions of the Schr\"odinger operator (\ref{defL}).
\end{corollary}

\subsection*{Acknowledgments}
I want to thank Iryna Egorova and Gerald Teschl for many discussions on this topic and the referee for valuable suggestions for improvements and hints with respect to the literature.

\end{document}